\documentclass[12pt]{article}
\usepackage{amssymb}
\usepackage{amsfonts}
\usepackage{amsmath}
\usepackage{amsbsy}
\newcommand{\be}{\begin{equation}}
\newcommand{\ee}{\end{equation}}
\newcommand{\bea}{\begin{eqnarray}}
\newcommand{\eea}{\end{eqnarray}}
\newcommand{\ba}{\begin{array}}
\newcommand{\ea}{\end{array}}

\newcommand{\bc}{\begin{center}}
\newcommand{\ec}{\end{center}}
\newcommand{\ben}{\begin{enumerate}}
\newcommand{\een}{\end{enumerate}}
\newcommand{\bfi}{\begin{figure}}
\newcommand{\efi}{\end{figure}}

\newcommand{\bq}{\begin{quote}}
\newcommand{\eq}{\end{quote}}
\newcommand{\bqu}{\begin{quotation}}
\newcommand{\equ}{\end{quotation}}
\newenvironment{emphit}{\begin{itemize}}{\end{itemize}}
\newcommand{\bemp}{\begin{emphit}}
\newcommand{\eemp}{\end{emphit}}

\newcommand{\bt}{\begin{tabular}}
\newcommand{\et}{\end{tabular}}

\newtheorem{myth}{Theorem}[section]

\newtheorem{mycor}{Corollary}[section]

\newtheorem{mydef}{Definition}[section]
\newtheorem{myrem}{Remark}[section]

\begin{document}
\date{}
\title{An Holomorphic Study Of Smarandache Automorphic and Cross Inverse Property Loops
}
\author{T\`em\'it\d{\'o}p\d{\'e} Gb\d{\'o}l\'ah\`an Ja\'iy\'e\d ol\'a\thanks{All correspondence to be addressed to this author.} \\
Department of Mathematics,\\
Obafemi Awolowo University, Ile Ife, Nigeria.\\
jaiyeolatemitope@yahoo.com, tjayeola@oauife.edu.ng}
 \maketitle
\begin{abstract}
By studying the holomorphic structure of automorphic inverse
property quasigroups and loops[AIPQ and (AIPL)] and cross inverse
property quasigroups and loops[CIPQ and (CIPL)], it is established
that the holomorph of a loop is a  Smarandache; AIPL, CIPL, K-loop,
Bruck-loop or Kikkawa-loop if and only if its Smarandache
automorphism group is trivial and the loop is itself is a
Smarandache; AIPL, CIPL, K-loop, Bruck-loop or Kikkawa-loop.
\end{abstract}

\section{Introduction}
\subsection{Quasigroups And Loops}
Let $L$ be a non-empty set. Define a binary operation ($\cdot $) on
$L$ : If $x\cdot y\in L$ for all $x, y\in L$, $(L, \cdot )$ is
called a groupoid. If the system of equations ;
\begin{displaymath}
a\cdot x=b\qquad\textrm{and}\qquad y\cdot a=b
\end{displaymath}
have unique solutions for $x$ and $y$ respectively, then $(L, \cdot
)$ is called a quasigroup. For each $x\in L$, the elements $x^\rho
=xJ_\rho ,x^\lambda =xJ_\lambda\in L$ such that $xx^\rho=e^\rho$ and
$x^\lambda x=e^\lambda$ are called the right, left inverses of $x$
respectively. Now, if there exists a unique element $e\in L$ called
the identity element such that for all $x\in L$, $x\cdot e=e\cdot
x=x$, $(L, \cdot )$ is called a loop.  To every loop $(L,\cdot )$
with automorphism group $AUM(L,\cdot )$, there corresponds another
loop. Let the set $H=(L,\cdot )\times AUM(L,\cdot )$. If we define
'$\circ$' on $H$ such that $(\alpha, x)\circ (\beta,
y)=(\alpha\beta, x\beta\cdot y)$ for all $(\alpha, x),(\beta, y)\in
H$, then $H(L,\cdot )=(H,\circ)$ is a loop as shown in Bruck
\cite{phd82} and is called the Holomorph of $(L,\cdot )$.
\paragraph{}
A loop(quasigroup) is a weak inverse property loop
(quasigroup)[WIPL(WIPQ)] if and only if it obeys the identity
\begin{equation*}\label{eq:8}
x(yx)^\rho=y^\rho\qquad\textrm{or}\qquad(xy)^\lambda x=y^\lambda.
\end{equation*}
A loop(quasigroup) is a cross inverse property
loop(quasigroup)[CIPL(CIPQ)] if and only if it obeys the identity
\begin{equation*}\label{eq:8.1}
xy\cdot x^\rho =y\qquad\textrm{or}\qquad x\cdot yx^\rho
=y\qquad\textrm{or}\qquad x^\lambda\cdot
(yx)=y\qquad\textrm{or}\qquad x^\lambda y\cdot x=y.
\end{equation*}
A loop(quasigroup) is an automorphic inverse property
loop(quasigroup)[AIPL(AIPQ)] if and only if it obeys the identity
\begin{equation*}
(xy)^\rho=x^\rho y^\rho~or~(xy)^\lambda =x^\lambda y^\lambda
\end{equation*}
Consider $(G,\cdot )$ and $(H,\circ )$ being two distinct
groupoids(quasigroups, loops). Let $A,B$ and $C$ be three distinct
non-equal bijective mappings, that maps $G$ onto $H$. The triple
$\alpha =(A,B,C)$ is called an isotopism of $(G,\cdot )$ onto
$(H,\circ )$ if and only if
\begin{displaymath}
xA\circ yB=(x\cdot y)C~\forall~x,y\in G.
\end{displaymath}
The set $SYM(G, \cdot )=SYM(G)$ of all bijections in a groupoid
$(G,\cdot )$ forms a group called the permutation(symmetric) group
of the groupoid $(G,\cdot )$. If $(G,\cdot )=(H,\circ )$, then the
triple $\alpha =(A,B,C)$ of bijections on $(G,\cdot )$ is called an
autotopism of the groupoid(quasigroup, loop) $(G,\cdot )$. Such
triples form a group $AUT(G,\cdot )$ called the autotopism group of
$(G,\cdot )$. Furthermore, if $A=B=C$, then $A$ is called an
automorphism of the groupoid(quasigroup, loop) $(G,\cdot )$. Such
bijections form a group $AUM(G,\cdot )$ called the automorphism
group of $(G,\cdot )$.

The left nucleus of $L$ denoted by $N_\lambda (L, \cdot )=\{a\in L :
ax\cdot y=a\cdot xy~\forall~x, y\in L\}$. The right nucleus of $L$
denoted by $N_\rho (L, \cdot )=\{a\in L : y\cdot xa=yx\cdot
a~\forall~ x, y\in L\}$. The middle nucleus of $L$ denoted by $N_\mu
(L, \cdot )=\{a\in L : ya\cdot x=y\cdot ax~\forall~x, y\in L\}$. The
nucleus of $L$ denoted by $N(L, \cdot )=N_\lambda (L, \cdot )\cap
N_\rho (L, \cdot )\cap N_\mu (L, \cdot )$. The centrum of $L$
denoted by $C(L, \cdot )=\{a\in L : ax=xa~\forall~x\in L\}$. The
center of $L$ denoted by $Z(L, \cdot )=N(L, \cdot )\cap C(L, \cdot
)$.

As observed by Osborn \cite{phd89}, a loop is a WIPL and an AIPL if
and only if it is a CIPL. The past efforts of Artzy \cite{phd140,
phd30, phd158, phd193}, Belousov and Tzurkan \cite{phd192} and
recent studies of Keedwell \cite{phd176}, Keedwell and Shcherbacov
\cite{phd175, phd177, phd178} are of great significance in the study
of WIPLs, AIPLs, CIPQs and CIPLs, their generalizations(i.e
m-inverse loops and quasigroups, (r,s,t)-inverse quasigroups) and
applications to cryptography. For more on loops and their
properties, readers should check \cite{phd41},\cite{phd39},
\cite{phd49}, \cite{phd42}, \cite{phd75} and \cite{phd3}.

Interestingly, Adeniran \cite{phd79} and Robinson \cite{phd85},
Oyebo and Adeniran \cite{phd141}, Chiboka and Solarin \cite{phd80},
Bruck \cite{phd82}, Bruck and Paige \cite{phd40}, Robinson
\cite{phd7}, Huthnance \cite{phd44} and Adeniran \cite{phd79} have
respectively studied the holomorphs of Bol loops, central loops,
conjugacy closed loops, inverse property loops, A-loops, extra
loops, weak inverse property loops, Osborn loops and Bruck loops.
Huthnance \cite{phd44} showed that if $(L,\cdot )$ is a loop with
holomorph $(H,\circ)$, $(L,\cdot )$ is a WIPL if and only if
$(H,\circ)$ is a WIPL. The holomorphs of an AIPL and a CIPL are yet
to be studied.

For the definitions of inverse property loop (IPL), Bol loop and
A-loop readers can check earlier references on loop theory.

Here ; a K-loop is an A-loop with the AIP, a Bruck loop is a Bol
loop with the AIP and a Kikkawa loop is an A-loop with the IP and
AIP.

\subsection{Smarandache Quasigroups And Loops}
The study of Smarandache loops was initiated by W. B. Vasantha
Kandasamy in 2002. In her book \cite{phd75}, she defined a
Smarandache loop (S-loop) as a loop with at least a subloop which
forms a subgroup under the binary operation of the loop. In
\cite{sma2}, the present author defined a Smarandache quasigroup
(S-quasigroup) to be a quasigroup with at least a non-trivial
associative subquasigroup called a Smarandache subsemigroup
(S-subsemigroup). Examples of Smarandache quasigroups are given in
Muktibodh \cite{muk2}. In her book, she introduced over 75
Smarandache concepts on loops. In her first paper \cite{phd83}, on
the study of Smarandache notions in algebraic structures, she
introduced Smarandache : left(right) alternative loops, Bol loops,
Moufang loops, and Bruck loops. But in \cite{sma1}, the present
author introduced Smarandache : inverse property loops (IPL), weak
inverse property loops (WIPL), G-loops, conjugacy closed loops
(CC-loop), central loops, extra loops, A-loops, K-loops, Bruck
loops, Kikkawa loops, Burn loops and homogeneous loops.

A loop is called a Smarandache A-loop(SAL) if it has at least a
non-trivial subloop that is a A-loop.

A loop is called a Smarandache K-loop(SKL) if it has at least a
non-trivial subloop that is a K-loop.

A loop is called a Smarandache Bruck-loop(SBRL) if it has at least a
non-trivial subloop that is a Bruck-loop.

A loop is called a Smarandache Kikkawa-loop(SKWL) if it has at least
a non-trivial subloop that is a Kikkawa-loop.

If $L$ is a S-groupoid with a S-subsemigroup $H$, then the set
$SSYM(L, \cdot )=SSYM(L)$ of all bijections $A$ in $L$ such that
$A~:~H\to H$ forms a group called the Smarandache
permutation(symmetric) group of the S-groupoid. In fact, $SSYM(L)\le
SYM(L)$.

The left Smarandache nucleus of $L$ denoted by $SN_\lambda (L, \cdot
)=N_\lambda (L, \cdot )\cap H$. The right Smarandache nucleus of $L$
denoted by $SN_\rho (L, \cdot )=N_\rho (L, \cdot )\cap H$. The
middle Smarandache nucleus of $L$ denoted by $SN_\mu (L, \cdot
)=N_\mu (L, \cdot )\cap H$. The Smarandache nucleus of $L$ denoted
by $SN(L, \cdot )=N(L, \cdot )\cap H$. The Smarandache centrum of
$L$ denoted by $SC(L, \cdot )=C(L, \cdot )\cap H$. The Smarandache
center of $L$ denoted by $SZ(L, \cdot )=Z(L, \cdot )\cap H$.

\begin{mydef}\label{1:1}
Let $(L, \cdot )$ and $(G, \circ )$ be two distinct groupoids that
are isotopic under a triple $(U, V, W)$. Now, if $(L, \cdot )$ and
$(G, \circ )$ are S-groupoids with S-subsemigroups $L'$ and $G'$
respectively such that $A~:~L'\to G'$, where $A\in\{U,V,W\}$, then
the isotopism $(U, V, W) : (L, \cdot )\rightarrow (G, \circ )$ is
called a Smarandache isotopism(S-isotopism).

Thus, if $U=V=W$, then $U$ is called a Smarandache isomorphism,
hence we write $(L, \cdot )\succsim (G, \circ )$.

But if $(L, \cdot )=(G, \circ )$, then the autotopism $(U, V, W)$ is
called a Smarandache autotopism (S-autotopism) and they form a group
$SAUT(L,\cdot )$ which will be called the Smarandache autotopism
group of $(L, \cdot )$. Observe that $SAUT(L,\cdot )\le AUT(L,\cdot
)$. Furthermore, if $U=V=W$, then $U$ is called a Smarandache
automorphism of $(L,\cdot )$. Such Smarandache permutations form a
group $SAUM(L,\cdot )$ called the Smarandache automorphism
group(SAG) of $(L,\cdot )$.
\end{mydef}
Let $L$ be a S-quasigroup with a S-subgroup $G$. Now, set
$H_S=(G,\cdot )\times SAUM(L,\cdot )$. If we define '$\circ$' on
$H_S$ such that $(\alpha, x)\circ (\beta, y)=(\alpha\beta,
x\beta\cdot y)$ for all $(\alpha, x),(\beta, y)\in H_S$, then
$H_S(L,\cdot )=(H_S,\circ)$ is a quasigroup.

If in $L$, $s^\lambda\cdot s\alpha\in SN(L)$ or $s\alpha\cdot
s^\rho\in SN(L)~\forall s\in G$ and $\alpha\in SAUM(L,\cdot  )$,
$(H_S,\circ)$ is called a Smarandache Nuclear-holomorph of $L$, if
$s^\lambda\cdot s\alpha\in SC(L)$ or $s\alpha\cdot s^\rho\in
SC(L)~\forall~s\in G$ and $\alpha\in SAUM(L,\cdot  )$, $(H_S,\circ)$
is called a Smarandache Centrum-holomorph of $L$ hence a Smarandache
Central-holomorph if $s^\lambda\cdot s\alpha\in SZ(L)$ or
$s\alpha\cdot s^\rho\in SZ(L)~\forall~s\in G$ and $\alpha\in
SAUM(L,\cdot )$.

The aim of the present study is to investigate the holomorphic
structure of Smarandache AIPLs and CIPLs(SCIPLs and SAIPLs) and use
the results to draw conclusions for Smarandache K-loops(SKLs),
Smarandache Bruck-loops(SBRLs) and Smarandache Kikkawa-loops
(SKWLs). This is done as follows.
\begin{enumerate}
\item The holomorphic structure of AIPQs(AIPLs) and CIPQs(CIPLs) are
investigated. Necessary and sufficient conditions for the holomorph
of a quasigroup(loop) to be an AIPQ(AIPL) or CIPQ(CIPL) are
established. It is shown that if the holomorph of a quasigroup(loop)
is a AIPQ(AIPL) or CIPQ(CIPL), then the holomorph is isomorphic to
the quasigroup(loop). Hence, the holomorph of a quasigroup(loop) is
an AIPQ(AIPL) or CIPQ(CIPL) if and only if its automorphism group is
trivial and the quasigroup(loop) is a AIPQ(AIPL) or CIPQ(CIPL).
Furthermore, it is discovered that if the holomorph of a
quasigroup(loop) is a CIPQ(CIPL), then the quasigroup(loop) is a
flexible unipotent CIPQ(flexible CIPL of exponent $2$).
\item The holomorph of a loop is shown to be a SAIPL, SCIPL, SKL, SBRL or SKWL
respectively if and only its SAG is trivial and the loop is a SAIPL,
SCIPL, SKL, SBRL, SKWL respectively.
\end{enumerate}

\section{Main Results}
\begin{myth}\label{3:1}
Let $(L,\cdot )$ be a quasigroup(loop) with holomorph $H(L)$. $H(L)$
is an AIPQ(AIPL) if and only if
\begin{enumerate}
\item $AUM(L)$ is an abelian group,
\item $(\beta^{-1}, \alpha ,I)\in AUT(L)~\forall~\alpha ,\beta\in AUM(L)$ and
\item $L$ is a AIPQ(AIPL).
\end{enumerate}
\end{myth}
{\bf Proof}\\
A quasigroup(loop) is an automorphic inverse property loop(AIPL) if
and only if it obeys the identity

Using either of the definitions of an AIPQ(AIPL), it can be shown
that $H(L)$ is a AIPQ(AIPL) if and only if $AUM(L)$ is an abelian
group and $(\beta^{-1}J_\rho, \alpha J_\rho ,J_\rho)\in
AUT(L)~\forall~\alpha ,\beta\in AUM(L)$. $L$ is isomorphic to a
subquasigroup(subloop) of $H(L)$, so $L$ is a AIPQ(AIPL) which
implies $(J_\rho, J_\rho ,J_\rho)\in AUT(L)$. So, $(\beta^{-1},
\alpha ,I)\in AUT(L)~\forall~\alpha ,\beta\in AUM(L)$.

\begin{mycor}\label{3:2}
Let $(L,\cdot )$ be a quasigroup(loop) with holomorph $H(L)$. $H(L)$
is a CIPQ(CIPL) if and only if
\begin{enumerate}
\item $AUM(L)$ is an abelian group,
\item $(\beta^{-1}, \alpha ,I)\in AUT(L)~\forall~\alpha ,\beta\in AUM(L)$ and
\item $L$ is a CIPQ(CIPL).
\end{enumerate}
\end{mycor}
{\bf Proof}\\
A quasigroup(loop) is a CIPQ(CIPL) if and only if it is a WIPQ(WIPL)
and an AIPQ(AIPL). $L$ is a WIPQ(WIPL) if and only if $H(L)$ is a
WIPQ(WIPL).

If $H(L)$ is a CIPQ(CIPL), then $H(L)$ is both a WIPQ(WIPL) and a
AIPQ(AIPL) which implies 1., 2., and 3. of Theorem~\ref{3:1}. Hence,
$L$ is a CIPQ(CIPL). The converse follows by just doing the reverse.

\begin{mycor}\label{3:3}
Let $(L,\cdot )$ be a quasigroup(loop) with holomorph $H(L)$. If
$H(L)$ is an AIPQ(AIPL) or CIPQ(CIPL), then $H(L)\cong L$.
\end{mycor}
{\bf Proof}\\
By 2. of Theorem~\ref{3:1}, $(\beta^{-1}, \alpha ,I)\in
AUT(L)~\forall~\alpha ,\beta\in AUM(L)$ implies $x\beta^{-1}\cdot
y\alpha =x\cdot y$ which means $\alpha =\beta =I$ by substituting
$x=e$ and $y=e$. Thus, $AUM(L)=\{I\}$ and so $H(L)\cong L$.

\begin{myth}\label{3:3.1}
The holomorph of a quasigroup(loop) $L$ is a AIPQ(AIPL) or
CIPQ(CIPL) if and only if $AUM(L)=\{I\}$ and $L$ is a AIPQ(AIPL) or
CIPQ(CIPL).
\end{myth}
{\bf Proof}\\
This is established using Theorem~\ref{3:1}, Corollary~\ref{3:2} and
Corollary~\ref{3:3}.

\begin{myth}\label{3:4}
Let $(L,\cdot )$ be a quasigroups(loop) with holomorph $H(L)$.
$H(L)$ is a CIPQ(CIPL) if and only if $AUM(L)$ is an abelian group
and any of the following is true for all $x,y\in L$ and $\alpha
,\beta\in AUM(L)$:
\begin{enumerate}
\item $(x\beta\cdot y)x^\rho=y\alpha$.
\item $x\beta\cdot yx^\rho=y\alpha$.
\item $(x^\lambda\alpha^{-1}\beta\alpha\cdot y\alpha )\cdot x=y$.
\item $x^\lambda\alpha^{-1}\beta\alpha\cdot (y\alpha \cdot x)=y$.
\end{enumerate}
\end{myth}
{\bf Proof}\\
This is achieved by simply using the four equivalent identities that
define a CIPQ(CIPL):

\begin{mycor}\label{3:5}
Let $(L,\cdot )$ be a quasigroups(loop) with holomorph $H(L)$. If
$H(L)$ is a CIPQ(CIPL) then, the following are equivalent to each
other
\begin{enumerate}
\item $(\beta^{-1}J_\rho, \alpha J_\rho ,J_\rho)\in AUT(L)~\forall~\alpha
,\beta\in AUM(L)$.
\item $(\beta^{-1}J_\lambda, \alpha J_\lambda ,J_\lambda )\in AUT(L)~\forall~\alpha
,\beta\in AUM(L)$.
\item $(x\beta\cdot y)x^\rho=y\alpha$.
\item $x\beta\cdot yx^\rho=y\alpha$.
\item $(x^\lambda\alpha^{-1}\beta\alpha\cdot y\alpha )\cdot x=y$.
\item $x^\lambda\alpha^{-1}\beta\alpha\cdot (y\alpha \cdot x)=y$.
\end{enumerate}
Hence,
\begin{displaymath}
(\beta, \alpha ,I),(\alpha, \beta ,I),(\beta,I, \alpha ),
(I,\alpha,\beta)\in AUT(L)~\forall~\alpha ,\beta\in AUM(L).
\end{displaymath}
\end{mycor}
{\bf Proof}\\
The equivalence of the six conditions follows from Theorem~\ref{3:4}
and the proof of Theorem~\ref{3:1}. The last part is simple.

\begin{mycor}\label{3:6}
Let $(L,\cdot )$ be a quasigroup(loop) with holomorph $H(L)$. If
$H(L)$ is a CIPQ(CIPL) then, $L$ is a flexible unipotent
CIPQ(flexible CIPL of exponent $2$).
\end{mycor}
{\bf Proof}\\
It is observed that $J_\rho =J_\lambda =I$. Hence, the conclusion
follows.

\begin{myrem}
The holomorphic structure of loops such as extra loop, Bol-loop,
C-loop, CC-loop and A-loop have been found to be characterized by
some special types of automorphisms such as
\begin{enumerate}
\item Nuclear automorphism(in the case of Bol-,CC- and extra loops),
\item central automorphism(in the case of central and A-loops).
\end{enumerate}
By Theorem~\ref{3:1} and Corollary~\ref{3:2}, the holomorphic
structure of AIPLs and CIPLs is characterized by commutative
automorphisms.
\end{myrem}

\begin{myth}\label{3:3.2}
The holomorph $H(L)$ of a quasigroup(loop) $L$ is a Smarandache
AIPQ(AIPL) or CIPQ(CIPL) if and only if $SAUM(L)=\{I\}$ and $L$ is a
Smarandache AIPQ(AIPL) or CIPQ(CIPL).
\end{myth}
{\bf Proof}\\
Let $L$ be a quasigroup with holomorph $H(L)$. If $H(L)$ is a
SAIPQ(SCIPQ), then there exists a S-subquasigroup $H_S(L)\subset
H(L)$ such that $H_S(L)$ is a AIPQ(CIPQ). Let $H_S(L)=G\times
SAUM(L)$ where $G$ is the S-subquasigroup of $L$. From
Theorem~\ref{3:3.1}, it can be seen that $H_S(L)$ is a AIPQ(CIPQ) if
and only if $SAUM(L)=\{I\}$ and $G$ is a AIPQ(CIPQ). So the
conclusion follows.

\begin{mycor}\label{}
The holomorph $H(L)$ of a loop $L$ is a SKL or SBRL or SKWL if and
only if $SAUM(L)=\{I\}$ and $L$ is a SKL or SBRL or SKWL.
\end{mycor}
{\bf Proof}\\
Let $L$ be a loop with holomorph $H(L)$. Consider the subloop
$H_S(L)$ of $H(L)$ such that $H_S(L)=G\times SAUM(L)$ where $G$ is
the subloop of $L$.
\begin{enumerate}
\item Recall that by [Theorem~5.3, \cite{phd40}], $H_S(L)$ is an
A-loop if and only if it is a Smarandache Central-holomorph of $L$
and $G$ is an A-loop. Combing this fact with Theorem~\ref{3:3.2}, it
can be concluded that: the holomorph $H(L)$ of a loop $L$ is a SKL
if and only if $SAUM(L)=\{I\}$ and $L$ is a SKL.
\item Recall that by \cite{phd85} and \cite{phd79}, $H_S(L)$ is a
Bol loop if and only if it is a Smarandache Nuclear-holomorph of $L$
and $G$ is a Bol-loop. Combing this fact with Theorem~\ref{3:3.2},
it can be concluded that: the holomorph $H(L)$ of a loop $L$ is a
SBRL if and only if $SAUM(L)=\{I\}$ and $L$ is a SBRL.
\item Following the first reason in 1., and using Theorem~\ref{3:3.2},
it can be concluded that: the holomorph $H(L)$ of a loop $L$ is a
SKWL if and only if $SAUM(L)=\{I\}$ and $L$ is a SKWL.
\end{enumerate}

\end{document}